\documentclass{amsart}[12pt]
\usepackage{amsmath,mathrsfs,amssymb}
\newtheorem{pro}{Proposition}
\newtheorem{nota}{Remark}

\title{On positive classical and free stable laws}
\begin{document}
\maketitle
\centerline{N. Demni\footnote{IRMAR, Rennes 1 University, France. Email: nizar.demni@univ-rennes1.fr
}}
\begin{abstract}
We derive the representative Bernstein measure of the density of $(X_{\alpha})^{-\alpha/(1-\alpha)}, 0 < \alpha < 1$, where $X_{\alpha}$ is a positive stable random variable, as a Fox-H function. When $1-\alpha = 1/j$ for some integer $j \geq 2$, the Fox H-function reduces to a Meijer G-function so that the Kanter's random variable (see below) is closely related to a product of $(j-1)$ independent Beta random variables. When $\alpha$ tends to $0$, the Bernstein measure becomes degenerate thereby agrees with Cressie's result for the asymptotic behaviour of stable distributions for small values of $\alpha$. Coming to free probability, our result makes more explicit that of Biane on the density of its free analog. The paper is closed with analytic arguments explaining the occurence of the Kanter's random variable in both the classical and the free settings.    				
\end{abstract}

\section{Motivation: Kanter's random variable}
The study of stable random variables has known a considerable growth during approximately the last five decades, since they are shown to be an efficient model for various phenomena occuring in biology, quantum physics, market finance (\cite{Zolo1}). Unfortunately, their densities has no closed formulas except in few cases. Consider for instance, the stable random variable $X_{\alpha}$ of index $0 < \alpha < 1$ and of asymmetry parameter equal one, known also as the positive stable variable since it is supported by $(0,\infty)$ (\cite{Ibr}, p.50). Its Laplace transform is given by $e^{-t^{\alpha}}, t > 0$ and its density admits the following expansion ((2.41) p.54 in \cite{Ibr}, (14.31) p.88 in \cite{Sato}, see also \cite{Zolo2}):
\begin{equation*} 
-\frac{1}{\pi x} \sum_{k=1}^{\infty} \frac{(-1)^k}{k!} \sin(k\pi \alpha) \Gamma(1+k\alpha) \frac{1}{x^{k\alpha}}, \quad x > 0.
\end{equation*}
The positivity of this series is not trivial and follows for instance from the integral representation displayed p.74 in \cite{Zolo1}. Amazing and not trivial as well is the complete monotonicity of the density of $(X_{\alpha})^{-r}, r \geq \alpha/(1-\alpha)$ (in particular its infinite divisibility, \cite{Sato}) which follows from the following representation (\cite{PSS}): there exists a standard exponential random variable $L$ and an infinitely divisible random variable $V_{\alpha,r}$, both variables being independent and 
\begin{equation*}
(X_{\alpha})^{r} \overset{d}{=} \frac{e^{V_{\alpha,r}}}{L}.
\end{equation*}
Recall that the proof of this fact relies on a suitable integral representation of the Gamma function and that the random variable $V_{\alpha,r}$ is characterized by its Laplace transform. For $r = \alpha/(1-\alpha)$, the above representation is more precise and traces back implicitely to Chernin and Ibragimov's paper \cite{Chernin} (see also Ch.II in \cite{Ibr}). In fact (\cite{Kanter} p.703) 
\begin{equation*}
(X_{\alpha})^{\alpha/(1-\alpha)} \overset{d}{=} \frac{a_{\alpha}(U)}{L},
\end{equation*}
where $U$ is a uniform random variable on $(0,\pi)$ and 
\begin{equation*}
a_{\alpha}(u) = [b_{\alpha}(u)]^{1/(1-\alpha)}, b(u) = [c_{\alpha}(u)]^{\alpha}[c_{1-\alpha}(u)]^{(1-\alpha)}, c_{\alpha}(u):= \frac{\sin(\alpha u)}{\sin(u)}.
\end{equation*}
We shall refer to $a_{\alpha}(U)$ as the Kanter's random variable. From the injectivity of the Laplace transform, one gets 
\begin{equation*}
V_{\alpha,1-\alpha} \overset{d}{=} \log[a_{\alpha}(U)].
\end{equation*}
Mysteriously, the function $a_{\alpha}$ is intimately related to the density of the free stable law of index $0 < \alpha < 1$ and asymmetry parameter $\rho=1$ derived p.1053 in \cite{Berco}, that we shall call shortly the positive free stable law. More precisely, if $U$ is defined on some probability space $(\Omega,\mathscr{F},\mathbb{P})$, then 
\begin{equation*}
\mathbb{P}(a_{\alpha}(U) \leq u) = \frac{1}{\pi} a_{\alpha}^{-1}(u), \quad u \geq a_{\alpha}(0) =  (1-\alpha)\alpha^{\alpha/(1-\alpha)}
\end{equation*}
since $a_{\alpha}$ is strictly increasing (\cite{Kanter} p.704), while the density of the positive free stable distribution takes the form (\cite{Berco}): 
\begin{equation}\label{Biane}
\frac{1}{\pi x} \frac{\sin(a_{1-\alpha}^{-1}(x))\sin(\alpha a_{1-\alpha}^{-1}(x))}{\sin((1-\alpha) a_{1-\alpha}^{-1}(x))}, 
\end{equation}
for $x \geq a_{1-\alpha}(0)$. That is why we thought it is reasonable to start with seeking an expression for the density of $a_{\alpha}(U)$ as a special function and more generally for that of $e^{V_{\alpha,r}}$. The easiest way to proceed is to invert the Mellin transform of the density of $e^{-V_{\alpha,r}}$ given by (\cite{PSS} p.292): 
\begin{equation*}
\mathbb{E}(e^{-sV_{\alpha,r}}) = \frac{\Gamma(rs/\alpha +1)}{\Gamma(s+1)\Gamma(sr +1)}, \quad \Re(s) > -\alpha/r \geq -(1-\alpha).
\end{equation*}
However, the Mellin's inversion formula 
\begin{equation}\label{Mellin}
\frac{1}{2i\pi y}\int_{-i\infty}^{+i\infty} \frac{\Gamma(rs/\alpha +1)}{\Gamma(s+1)\Gamma(sr +1)} y^{-s} ds,
\end{equation}
applies for any $r > \alpha/(1-\alpha)$ (we take $\Re(s) = 0$ as a path of integration) and does not when $r=\alpha/(1-\alpha)$. The latter fact is seen from the following estimate of $|\Gamma(a+ib)|, a,b \in \mathbb{R}, |b| \rightarrow \infty$ displayed in formula 1.24. p.4 from \cite{HMS}: 
\begin{equation}\label{EST}
|\Gamma(a+ib)| \sim \sqrt{2\pi} |b|^{a-1/2}e^{-a-\pi |b|/2}, \quad |b| \rightarrow \infty
\end{equation} 
so that the absolute value of the integrand in \eqref{Mellin} is equivalent to $|s|^{-1/2}$ (up to a constant). Note however that when $r > \alpha/(1-\alpha)$, the integral displayed in \eqref{Mellin} is nothing else but a so-called Fox H-function (\cite{HMS} p.3, see below): 
\begin{equation*}
\frac{1}{2i\pi y}\int_{-i\infty}^{+i\infty} \frac{\Gamma(rs/\alpha +1)}{\Gamma(s+1)\Gamma(sr +1)} y^{-s} ds = \frac{1}{y} H_{2,1}^{1,0}\left[y|_{(1,r/\alpha)}^{(1,1),(1,r)}\right]
\end{equation*}
so that the density of $e^{V_{\alpha,r}}$ reads 
\begin{equation*}
\frac{1}{y} H_{2,1}^{1,0}\left[\frac{1}{y}|_{(1,r/\alpha)}^{(1,1),(1,r)}\right], \quad y > 0, \, r > \frac{\alpha}{1-\alpha}. 
\end{equation*}
The main result of this paper states then that in the pathological case $r=\alpha/(1-\alpha)$, the density of $a_{\alpha}(U)$ is a H-function too represented through a path integral resembling to \eqref{Mellin} but with a different contour (Hankel path). When $\alpha$ approaches zero, this density tends to zero pointwisely, a fact that agrees with Cressie's result implying that the Kanter's random variable becomes degenerate in the limit $\alpha \rightarrow 0^+$. When $1-\alpha = 1/j$ for some integer $j \geq 2$, the derived H-function reduces to the so-called Meijer G-function so that $a_{\alpha}(U)$ is closely related to a product of suitably chosen $j-1$ independent Beta random variables. The latter fact resembles to Williams respresentation of positive stable variables (\cite{Will}). Coming to free probability theory, the distribution function of $a_{\alpha}(U)$ is computed yielding an explicit expression for the density of a positive free stable law. Using the free analog of Zolotarev's duality (\cite{Berco}), a similar result holds for a free stable law with index $1 < \alpha < 2$ and asymmetry parameter $\rho = 0$. The paper is closed with supplying analytic arguments that explain the occurence of $a_{\alpha}$ in both probabilistic settings. It turns out that it is resumed in the inverse formula for the Fourier transform of the density of $X_{\alpha}$ together with deformations of paths of integration. Nevertheless, we think that a group theoretical argument exists but seems to be very hidden. For sake of completeness, some basic needed facts on the H-function are collected in the next section. A good reference to them is the monograph \cite{HMS}. 

\section{On the Fox H-function}
The Fox H-function is defined as a Mellin-Barnes integral: 
\begin{equation*}
H_{p,q}^{m,n}\left[z|_{(b_i,B_i)_{1 \leq i \leq q}}^{(a_i,A_i)_{1 \leq i \leq p}}\right] = \int_L \Theta(s)z^{-s} ds
\end{equation*}
where $1 \leq m \leq q, 0\leq n \leq p, a_i,b_i \in \mathbb{R}, A_i,B_i > 0$, 
\begin{equation*}
\Theta(s) = \frac{\prod_{i=1}^m \Gamma(b_i + B_is)\prod_{i=1}^n\Gamma(1-a_i - A_is)}{\prod_{i=m+1}^q \Gamma(1- b_i - B_is)\prod_{i=n+1}^p\Gamma(a_i + A_is)},
\end{equation*}
the principal determination of the power function is taken (though it is not necessary here), an empty product being equal one and $L$ is a suitable contour separating the poles of both products of the meromorphic (Gamma) functions in the numerator of $\Theta$ (\cite{HMS} p.3,4). The choice of $L$ and the domain of convergence of the Mellin-Barnes integral defining the H-function for each contour depend on the parameters $a_i,b_i,A_i,B_i$. Eight cases are discussed in \cite{HMS} p.4. and we will only make use of five of them in the sequel that we shall refer to whenever needed: 
\begin{itemize}
\item[i)] $q \geq 1$,
\begin{equation*}
\mu:= \sum_{i=1}^q B_i - \sum_{i=1}^p A_i > 0,
\end{equation*}
then the H-function exists in the punctured complex plane.
\item[ii)] $p \geq 1, \mu = 0$, the H-function exists in the domain $|z| > \beta$ where 
\begin{equation*}
\beta:= \prod_{i=1}^p A_i^{-A_i} \prod_{i=1}B_i^{B_i}
\end{equation*}
and $L = L_{+\infty}$ described in \cite{HMS} p.3. 
\item[iii)] Let $\Omega$ be the real number defined by 
\begin{equation*}
\Omega:= \sum_{i=1}^n A_i + \sum_{i=1}^m B_i -  \sum_{i=n+1}^p A_i - \sum_{i=m+1}^q B_i.
\end{equation*}
Then the Fox H-function exists for all $z \neq 0$ such that $|\arg(z)| < \pi \Omega/2, \Omega > 0$ and $L$ is the infinite semi-circle $\gamma-i\infty, \gamma + i \infty$ for a suitable real $\gamma$.
\item[iv)] When $q \geq 1, \mu=0$, the H-function exists in the open disc $|z| < \beta$ and $L= L_{-\infty}$ is described in \cite{HMS} p.3.
\item[v)] Let 
\begin{equation*}
\delta:= \sum_{i=1}^q b_i - \sum_{i=1}^pa_i + \frac{p-q}{2}. 
\end{equation*}
If $\mu = 0,\delta < -1$, then the H-function is well-defined for complex numbers lying on the circle $|z| = \beta$. 
\end{itemize}
The most crucial fact when looking at these definitions is that whenever more than one definition make sense, they lead to the same value of the integral. In fact, the integral is evaluated by means of Cauchy's residue Theorem applied to truncated contours together with asymptotic estimates at infinity, and the contour of integration can be deformed into another one provided that the deformation is admissible (\cite{Dieud}, \cite{Whitt}, see also CH.V. in \cite{Erd} for a similar discussion on Meijer's G-function which fits the Fox H-function when $A_i = B_j = 1, 1 \leq i \leq p, 1 \leq j \leq q$).  

\begin{nota}
In the sequel, the complex variable $z$ will take only strictly positive values thereby $s \mapsto z^{-s}$ is defined through the real-valued logarithm function. 
\end{nota}
\section{The density of $a_{\alpha}(U)$}
According to Kanter's representation, there exists a density $h_{\alpha}$  supported by $a_{\alpha}(0,\pi) = (a_{\alpha}(0),\infty)$  such that 
\begin{equation*}
-\frac{1-\alpha}{\alpha\pi x} \sum_{k=1}^{\infty} \frac{(-1)^k}{k!} \sin(k\pi \alpha) \Gamma(1+k\alpha) x^{k(1-\alpha)}  = \int_0^{\infty}e^{-xy} y h_{\alpha}(y)dy.
\end{equation*}
The main result of this paper is stated as: 
\begin{pro}
The density of the Kanter's variable is expressed as: 
\begin{align*}
h_{\alpha}(y) = \frac{1}{y}H^{1,0}_{2,1}\left[\frac{1}{y}|_{(1,1/(1-\alpha))}^{(1,\alpha/(1-\alpha)),(1,1)}\right] = \frac{2(1-\alpha)}{y} H^{1,0}_{2,1}\left[\frac{1}{y^{2(1-\alpha)}}|_{(1,2)}^{(1,2\alpha),(1,2(1-\alpha))}\right]
\end{align*}
for $y \in [a_{\alpha}(0), \infty[$. 
 \end{pro}
{\it Proof}: We shall start from the density of $(X_{\alpha})^{-\alpha}$ which has the expansion
\begin{equation*}
-\frac{1}{\alpha\pi x} \sum_{k=1}^{\infty} \frac{(-1)^k}{k!} \sin(k\pi \alpha) \Gamma(1+k\alpha) x^{k}.
\end{equation*}
Using the mirror formula satisfied by the Euler's Gamma function (\cite{Erd} p.3),
\begin{equation*}
\Gamma(1+ k\alpha)\Gamma(-k\alpha) = - \frac{\pi}{\sin(k\alpha \pi)},
\end{equation*} 
the density of $X_{\alpha}^{-\alpha}$ is transformed to
\begin{equation*}
\frac{1}{\alpha x} \sum_{k=1}^{\infty} \frac{(-1)^k}{k!} \frac{x^k}{\Gamma(-k\alpha)}
\end{equation*}
and the infinite summation may start from $k=0$ since $z=0$ is a pole of the Gamma function. According to A.6 p.218 in \cite{HMS}, this is expressed through the H-function as follows
\footnote{It is a special kind of what V. M. Zolotarev called incomplete hypergeometric function in \cite{Zolo2}.}: 
\begin{align*}
-\frac{1}{\alpha\pi x} \sum_{k=1}^{\infty} \frac{(-1)^k}{k!} \sin(k\pi \alpha) \Gamma(1+k\alpha) x^{k} & = \frac{2}{\alpha x}H_{1,1}^{1,0}\left[x^2|_{(0,2)}^{(0,2\alpha)}\right]
\end{align*}
for all $x > 0$ (note that the H-function in the RHS exists for all $z \neq 0$ in the complex plane by i) since  $q \geq 1, \mu = 2-2\alpha > 0$). Therefore,  the density of 
$(X_{\alpha})^{-\alpha/(1-\alpha)}$ reads
\begin{equation*}
\frac{2(1-\alpha)}{\alpha x}H_{1,1}^{1,0}\left[x^{2(1-\alpha)}|_{(0,2)}^{(0,2\alpha)}\right].
\end{equation*}
According to iii), one has 
\begin{equation*}
H_{1,1}^{1,0}\left[x^{2(1-\alpha)}|_{(0,2)}^{(0,2\alpha)}\right] = \frac{1}{2i\pi} \int_{\gamma -i\infty}^{\gamma+i\infty} \frac{\Gamma(2s)}{\Gamma(2\alpha s)}  \frac{ds}{x^{2(1-\alpha)s}} 
\end{equation*}
where $\Re(s) = \gamma > 0$. With the help of the Gamma integral 
\begin{align*}
\frac{2(1-\alpha)}{\alpha x^2}H_{1,1}^{1,0}\left[x^{2(1-\alpha)}|_{(1,2)}^{(1,2\alpha)}\right] &= \frac{2(1-\alpha)}{2\alpha i\pi} \int_{\gamma-i\infty}^{\gamma+i\infty} \frac{\Gamma(2s)}{\Gamma(2\alpha s)\Gamma(2+2(1-\alpha)s)}  \int_0^{\infty}e^{-xy} y^{2(1-\alpha)s+1} dyds
\\& = \frac{2(1-\alpha)}{2 i\pi} \int_{\gamma-i\infty}^{\gamma+i\infty} \frac{\Gamma(2s+1)}{\Gamma(2\alpha s+1)\Gamma(2+2(1-\alpha)s)}  \int_0^{\infty}e^{-xy} y^{2(1-\alpha)s+1} dyds.
\end{align*}
In order to change the order of integration, one needs to check that 
\begin{equation*}
\int_{-\infty}^{+\infty} \frac{|\Gamma(2is +2\gamma+1)|}{|\Gamma(2i\alpha s+ 2\alpha \gamma + 1)\Gamma(2(1-\alpha)s+ 2(1-\alpha)\gamma + 2)|}  \int_0^{\infty}e^{-xy} y^{2(1-\alpha)\gamma} dyds < \infty.
 \end{equation*}
Performing an integration with respect to $y$, it amounts to check 
\begin{equation*}
\int_{-\infty}^{+\infty} \frac{|\Gamma(2is +2\gamma+1)|}{|\Gamma(2i\alpha s+ 2\alpha \gamma + 1)\Gamma(2i(1-\alpha)s+ 2(1-\alpha)\gamma + 2)|} ds < \infty
\end{equation*}
which follows again from the estimate of $|\Gamma(a+ib)|, |b| \rightarrow \infty$ (see \eqref{EST}). Thus, Fubini's Theorem yields
\begin{align}
\int_0^{\infty}e^{-xy}yh_{\alpha}(y)dy &= 2(1-\alpha) \int_0^{\infty}xe^{-xy}  \int_{\gamma-i\infty}^{\gamma+i\infty} \frac{\Gamma(2s+1)}{\Gamma(2\alpha s+1)\Gamma(2+2(1-\alpha)s)} 
y^{2(1-\alpha)s+1}ds dy \nonumber
\\& = 2(1-\alpha) \int_0^{\infty}xe^{-xy} yH^{1,0}_{2,1}\left[\frac{1}{y^{2(1-\alpha)}}|_{(1,2)}^{(1,2\alpha),(2,2(1-\alpha))}\right]dy. \label{E1}
\end{align}
Note that the H-function displayed in \eqref{E1} may be represented for $y > (1-\alpha)\alpha^{\alpha/(1-\alpha)}$ as a path integral through the loop $L_{-\infty}$ according to iv), it vanishes for $0 < y < (1-\alpha)\alpha^{\alpha/(1-\alpha)}$ since according to ii) the contour is a loop $L_{+\infty}$ that contains no pole of the meromorphic function $s \mapsto \Gamma(1+2s)$, and it takes a finite value at $z=(1-\alpha)\alpha^{\alpha/(1-\alpha)}$ regarding v) since $\mu=0, \delta= -3/2 < -1$ (in particular, it vanishes there by continuity). Hence, an integration by parts shows that the density of the Kanter's random variable is 
\begin{align*}
h_{\alpha}(y) &= \frac{2(1-\alpha)}{y} \frac{d}{dy}\left\{yH^{1,0}_{2,1}\left[\frac{1}{y^{2(1-\alpha)}}|_{(1,2)}^{(1,2\alpha),(2,2(1-\alpha))}\right]\right\}
\\& = \frac{2(1-\alpha)}{y} \frac{d}{dy}\int_{L_{-\infty}} \frac{\Gamma(2s+1)}{\Gamma(2\alpha s+1)\Gamma(2+2(1-\alpha)s)} y^{2(1-\alpha)s+1}ds
\\& = \frac{2(1-\alpha)}{y} \int_{L_{-\infty}} \frac{\Gamma(2s+1)}{\Gamma(2\alpha s+1)\Gamma(2(1-\alpha)s+1)} y^{2(1-\alpha)s}ds
\\& = \frac{2(1-\alpha)}{y} H^{1,0}_{2,1}\left[\frac{1}{y^{2(1-\alpha)}}|_{(1,2)}^{(1,2\alpha),(1,2(1-\alpha))}\right]
\\& =  \frac{1}{y} \int_{L'_{-\infty}} \frac{\Gamma(s/(1-\alpha)+1)}{\Gamma(\alpha s/(1-\alpha)+1)\Gamma(s+1)} y^{s}ds = 
 \frac{1}{y}H^{1,0}_{2,1}\left[\frac{1}{y}|_{(1,1/(1-\alpha))}^{(1,\alpha/(1-\alpha)),(1,1)}\right],
\end{align*}
where $L'_{-\infty}$ is the image of $L_{-\infty}$ under the map $s \mapsto 2(1-\alpha)s$. The density $h_{\alpha}$ is well defined for $y > (1-\alpha)\alpha^{\alpha/(1-\alpha)}$ by iv) and vanishes for $0 < y  < (1-\alpha)\alpha^{\alpha/(1-\alpha)}$ according to ii) (see also 1.33 p.6 in \cite{HMS}). $\hfill \blacksquare$

\section{Special indices} 
\subsection{The $1/2$-stable case} 
When $\alpha=1/2$, the density of the positive stable random variable reads (\cite{Zolo1} p.66)
\begin{equation*}
\frac{1}{2\sqrt{\pi}} x^{-3/2} e^{-1/(4x)} {\bf 1}_{\{x > 0\}}.
\end{equation*}
According to Kanter's representation, the image of the above density under the map $x \mapsto 1/x$, that is 
\begin{equation*}
\frac{1}{2\sqrt{\pi x}}e^{-x/4} {\bf 1}_{\{x > 0\}} = \frac{1}{2\pi}\sqrt{\frac{\pi}{x}}{\bf 1}_{\{x > 0\}} e^{-x/4}, 
\end{equation*}
should be the Laplace transform of $y \mapsto yh_{1/2}(y)$ which takes the form (after the use of the mirror formula): 
\begin{equation*}
yh_{1/2}(y) = \frac{1}{2\pi^2} \sum_{k=0}^{\infty}\sin^2[(k +1)\pi/2] \Gamma^2[(k+1)/2] \frac{(-1)^k}{k!} \frac{1}{y^{(k+1)/2}}, \, y > 1/4.
\end{equation*}
Note that $x \mapsto x^{-1/2}$ is completely monotone since $x \mapsto x^{1/2}$ is a Bernstein function, and that the representing measure is easily computed using the Gamma integral as 
\begin{equation*}
\frac{dy}{\sqrt{\pi y}} {\bf 1}_{\{ y > 0\}}. 
\end{equation*}
Now, since $\sin (n\pi) = 0$ for any integer $n$, then 
\begin{align*}
yh_{1/2}(y) &= \frac{1}{2\pi^2} \sum_{k=0}^{\infty} \frac{\sin^2[(2k +1)\pi/2] \Gamma^2[(2k+1)/2]}{(2k!)}\frac{1}{y^{(2k+1)/2}} 
\\& = \frac{1}{2\pi^2 y^{1/2}} \sum_{k=0}^{\infty} \frac{ \Gamma^2(k+1/2)}{\Gamma(2k+1)} \frac{1}{y^{k}}
\end{align*}
for $y > 1/4$. Using Legendre's duplication formula (\cite{Erd} p.5)
\begin{equation*}
\sqrt{\pi} \Gamma(2k+1) = 2^{2k} \Gamma(k+1/2) \Gamma(k+1),
\end{equation*}
one gets
\begin{align*}
yh_{1/2}(y)  &= \frac{1}{2\pi y^{1/2}}\sum_{k=0}^{\infty} \frac{(1/2)_k}{k!} \frac{1}{(4y)^{k}} = \frac{1}{2 \pi y^{1/2}} \frac{1}{\sqrt{1-1/(4y))}} \\
 & =  \frac{1}{2\pi} \frac{1}{\sqrt{y-1/4}}, \quad y > 1/4,
 \end{align*}
which is nothing but 
\begin{equation*}
\frac{1}{2\pi}\frac{1}{\sqrt{y}} {\bf 1}_{\{ y > 0\}} \star \delta_{1/4}(dy)
\end{equation*}
where $\star$ is the classical convolution of positive measures.

\subsection{Relation to product of independent Beta variables}
Let 
\begin{equation*}
1-\alpha = \frac{1}{j}, \quad \alpha = \frac{j}{j-1}
\end{equation*}
for some integer $j \geq 2$, then 
\begin{equation*}
\mathbb{E}\left[\frac{1}{[a_{1-1/j}(U)]^s}\right] = \frac{\Gamma(js +1)}{\Gamma(s+1)\Gamma((j-1)s +1)} \quad \Re(s) >  -1/j.
\end{equation*}
The multiplication Theorem (generalization of Legendre's duplication formula, \cite{Erd} p.4)
\begin{equation*}
\Gamma(js+1) = \frac{1}{(\sqrt{2\pi})^{j-1}} j^{js+1/2}\prod_{k=1}^{j} \Gamma\left(s + \frac{k}{j}\right)
\end{equation*}
yields
\begin{equation*}
\mathbb{E}\left[\frac{1}{[a_{1-1/j}(U)]^s}\right] = \frac{1}{\sqrt{2\pi}} \sqrt{\frac{j}{j-1}} \frac{j^{js}}{(j-1)^{(j-1)s}}  \prod_{k=1}^{j-1} \frac{\Gamma\left(s + k/j\right)}{\Gamma(s+k/(j-1))}. 
\end{equation*}

Now, let $\beta_1, \beta_2, \cdots, \beta_{j-1}$ be independent Beta random variables of the first kind such that the density of $\beta_k$ reads:
\begin{equation*}
\frac{1}{B(k/j+1, k/(j(j-1)))} t^{k/j}(1-t)^{k/(j-1) - k/j - 1}{\bf 1}_{(0,1)}(t) 
\end{equation*}
where $B(\cdot,\cdot)$ stands for the Beta function. According to formula (4.11) p.122 in \cite{HMS}, 
\begin{align*}
\mathbb{E}[(\beta_1 \cdots \beta_{j-1})^{s-1}] &= \prod_{k=1}^{j-1}\frac{\Gamma(k/(j-1)+1)}{\Gamma(k/j+1)} \prod_{k=1}^{j-1}\frac{\Gamma(s+ k/j)}{\Gamma(s+ k/(j-1))}
\\& = \frac{1}{\sqrt{2\pi}} \sqrt{\frac{j-1}{j}} \frac{j^{j}}{(j-1)^{(j-1)}} \frac{\Gamma(j-1)}{\Gamma(j)} \prod_{k=1}^{j-1}\frac{\Gamma(s+ k/j)}{\Gamma(s+ k/(j-1))}
\\& = \frac{1}{\sqrt{2\pi}} \sqrt{\frac{1}{j(j-1)}} \frac{j^{j}}{(j-1)^{(j-1)}} \prod_{k=1}^{j-1}\frac{\Gamma(s+ k/j)}{\Gamma(s+ k/(j-1))}
\end{align*}
again by the multiplication Theorem. It follows that 
\begin{equation*}
\mathbb{E}\left[\frac{1}{[a_{1-1/j}(U)]^s}\right] = j  \left[\frac{j^{j}}{(j-1)^{(j-1)}}\right]^{s-1}\mathbb{E}[(\beta_1 \cdots \beta_{j-1})^{s-1}].
\end{equation*}
When $s \neq 0$, the last equality may be written as 
\begin{equation*}
\int_0^{\infty} u^{s-1} \mathbb{P}(a_{1-1/j}(U) < 1/u) du =  j  \mathbb{E}\left[\left(\frac{j^{j}}{(j-1)^{(j-1)}}V\beta_1 \cdots \beta_{j-1}\right)^{s-1}\right]
\end{equation*}
where $V$ is uniformly distributed in $(0,1)$ and is independent from the $\beta_k$s. 

\begin{nota}
The relation to product of independent Beta variables is not surprising since the H-function involved in the density of $h_{1-1/j}$ reduces, when $1-\alpha = 1/j$ to a Meijer's G-function (\cite{HMS} p.123). It somewhat matches with Williams result (\cite{Will}) if one keeps in mind the Gamma-Beta algebra. In fact, $1/X_{\alpha}$ is distributed, when $\alpha = 1/j, j \geq 2$ as a product of $j-1$ Gamma random variables $\gamma_k, 1 \leq k \leq j-1$ with densities 
\begin{equation*}
\frac{1}{\Gamma(k/j)}t^{k/j - 1} e^{-t} {\bf 1}_{\{t > 0\}}, \quad 1 \leq k \leq j-1.
\end{equation*}
Moreover, since 
\begin{equation*}
[a_{\alpha}(\theta)]^{(1-\alpha)/\alpha} = a_{1-\alpha}(\theta),
\end{equation*}
then Williams result implies
\begin{equation*}
j^{j} \gamma_1\cdots \gamma_k \overset{d}{=}  \frac{1}{X_{1/j}} \overset{d}{=} \frac{L^{j-1}}{a_{1-1/j}(U)}. 
\end{equation*} 
Note also that, for $\alpha = 1/(jd), j \geq 2d, d \geq 1$, a more general representation of positive stable laws is due to T. Simon (\cite{Simon}) and involves independent Gamma and Beta variables. 
\end{nota}  

\section{Free positive stable distribution}
Free stable distributions may be defined in a similar way as in the classical setting, when substituting the classical convolution of probability measures by Voiculescu's free convolution (see \cite{Berco}). As mentioned in the introduction, the inverse function (in the composition's sense) of $a_{1-\alpha}$ provides an explicit expression of the density of a stable law of index $0< \alpha < 1$ and asymmetry parameter $\rho=1$ in the free probability setting. More precisely, if $V$ is a uniform random variable on $(0,1)$ defined on some probability space $(\Omega,\mathscr{F},\mathbb{P})$, then 
\begin{equation*}
 \frac{1}{\pi}a_{1-\alpha}^{-1}(x) = \mathbb{P}(a_{1-\alpha}(\pi V) \leq x) 
\end{equation*}
for $x > \alpha(1-\alpha)^{(1-\alpha)/\alpha}$. When $\alpha = 1/2$, it simplifies to 
\begin{align*}
a_{1/2}^{-1}(x) &= \frac{1}{2}\int_{1/4}^x \frac{dy}{y\sqrt{y-1/4}} = \frac{1}{2}\int_0^{x-1/4} \frac{dy}{(y+1/4)\sqrt{y}}
\\& = 2\int_{0}^{x-1/4} \frac{2dy}{1+4y^2} = 2\arctan\left[2\sqrt{x-\frac{1}{4}}\right].
\end{align*}
Using the trigonometric identity 
\begin{align*}
\sin(\theta) = 2\frac{\tan(\theta/2)}{1+\tan^2(\theta/2)}
\end{align*}
one recovers the density given by Biane (\cite{Berco})
\begin{align*}
\frac{2}{\pi x} \sin[a_{1/2}^{-1}(x)] &= \frac{2}{\pi x} \frac{\tan(a_{1/2}^{-1}(x)/2)}{1+\tan^2(a_{1/2}^{-1}(x)/2)} = \frac{2}{\pi x} \frac{\sqrt{4x-1}}{1 + 4(x-1/4)} =  \frac{\sqrt{4x-1}}{2\pi x^2}.
\end{align*}
For general $\alpha \in (0,1)$, we use the integral representation 
\begin{equation*}
H^{1,0}_{2,1}\left[\frac{1}{y^{2\alpha}}|_{(1,2)}^{(1,2\alpha),(1,2(1-\alpha))}\right] = \frac{1}{2i\pi}\int_{L_{-\infty}}\frac{\Gamma(2s+1)}{\Gamma(2\alpha s+1)\Gamma(1+2(1-\alpha)s)} y^{2\alpha s} ds.
\end{equation*}
Now, assume for a while that Fubini's Theorem applies so that 
\begin{align*}
 \frac{1}{\pi}a_{1-\alpha}^{-1}(x) &= \frac{2\alpha}{2i\pi}\int_{L_{-\infty}}\frac{\Gamma(2s+1)}{\Gamma(2\alpha s+1)\Gamma(1+2(1-\alpha)s)} \int_{\alpha(1-\alpha)^{(1-\alpha)/\alpha}}^{x}
 y^{2\alpha s-1} dy ds \\&
= \frac{2}{2i\pi}\int_{L_{-\infty}}\frac{\Gamma(2s)}{\Gamma(2\alpha s+1)\Gamma(1+2(1-\alpha)s)} \left[x^{2\alpha s} - (\alpha^{2\alpha}(1-\alpha)^{2(1-\alpha)})^s\right]ds.
 \end{align*}
Since the poles of $s \mapsto \Gamma(2s)$ are $s= -k/2, k \in \mathbb{N}$ and simple, and since the pole $s=0$ lies outside $L_{-\infty}$, then Cauchy's Residue Theorem yields 
\begin{equation*}
\frac{1}{2i\pi}\int_{L_{-\infty}}\frac{\Gamma(2s)}{\Gamma(2\alpha s+1)\Gamma(1+2(1-\alpha)s)} x^{2\alpha s} ds = H_{2,1}^{1,0}\left[\frac{1}{x^{2\alpha}}|_{(0,2)}^{(1,2\alpha),(1,2(1-\alpha))}\right]  - \textrm{Res}(x,0)
\end{equation*}
where $\textrm{Res}(x,0)$ is the residue of the meromorphic function 
\begin{equation*}
s \mapsto \frac{\Gamma(2s)}{\Gamma(2\alpha s+1)\Gamma(1+2(1-\alpha)s)} x^{2\alpha s}
\end{equation*}
at $s=0$ for fixed parameter $x >  \alpha(1-\alpha)^{(1-\alpha)/\alpha}$. This residue is easily computed as 
\begin{equation*}
\lim_{s \rightarrow 0} \frac{s\Gamma(2s)}{\Gamma(2\alpha s+1)\Gamma(1+2(1-\alpha)s)} x^{2\alpha s} = \frac{1}{2}, \quad x > \alpha^{2\alpha}(1-\alpha)^{2(1-\alpha)}. 
\end{equation*}
Since the last H-function is defined on the whole real line, in particular at $\alpha(1-\alpha)^{(1-\alpha)/\alpha}$ by the virtue of v) since $\mu=0$ and $\delta < -3/2 < -1$, then 
\begin{align*}
\frac{1}{\pi}a_{1-\alpha}^{-1}(x) = 2H_{2,1}^{1,0}\left[\frac{1}{x^{2\alpha}}|_{(0,2)}^{(1,2\alpha),(1,2(1-\alpha))}\right] - 2H_{2,1}^{1,0}\left[\frac{1}{\alpha^{2\alpha}(1-\alpha)^{2(1-\alpha)}}|_{(0,2)}^{(1,2\alpha),(1,2(1-\alpha))}\right]
 \end{align*}
Regarding the continuity of the H-function at $\alpha(1-\alpha)^{(1-\alpha)/\alpha}$, it vanishes there (since it vanishes for $x < \alpha(1-\alpha)^{(1-\alpha)/\alpha}$) and one gets 
\begin{align*}
\frac{1}{\pi}a_{1-\alpha}^{-1}(x) = 2H_{2,1}^{1,0}\left[\frac{1}{x^{2\alpha}}|_{(0,2)}^{(1,2\alpha),(1,2(1-\alpha))}\right], x > \alpha(1-\alpha)^{(1-\alpha)/\alpha} 
\end{align*}
and zero otherwise. Finally, in order to fit in the density of the free positive stable law displayed in \eqref{Biane}, we can rewrite the last H-function as 
\begin{equation*}
\frac{1}{\alpha} H_{2,1}^{1,0}\left[\frac{1}{x}|_{(0,1/\alpha)}^{(1,1),(1,(1-\alpha)/\alpha)}\right] = \frac{1}{1-\alpha} H_{2,1}^{1,0}\left[\frac{1}{x^{\alpha/(1-\alpha)}}|_{(0,1/(1-\alpha))}^{(1,\alpha/(1-\alpha),(1,1)}\right]. 
\end{equation*}
Note that the duality relation given in \cite{Berco} allows to derive the density of a free stable distribution of index $1 < \alpha < 2$ and asymmetry coefficient $\rho=0$ (the case of an asymmetry parameter $\rho = 1$ is obtained by the simple variable change $x \mapsto -x$). 
Coming back to the validity of Fubini's Theorem, one has to prove, after integrating with respect to the variable $y$, the convergence of 
\begin{equation*}
\int_{L_{-\infty}}\frac{|\Gamma(2s+1)|}{|\Gamma(2\alpha s+1)\Gamma(1+2(1-\alpha)s)|}  \left[x^{2\alpha \Re(s)} - (\alpha^{2\alpha}(1-\alpha)^{2(1-\alpha)})^{\Re(s)}\right] \frac{|ds|}{|\Re(s)|}. 
\end{equation*}
This is easily cheked from the estimate of $\Gamma(a+ib)$ for large $a$ (formula 1.23 p.4. in \cite{HMS}):
\begin{align*}
|\Gamma(a+ib)|& \sim \sqrt{2\pi} |a|^{a-1/2}e^{-a- a(1-\textrm{sgn}(b))/2}, \quad |a| \rightarrow \infty
\end{align*} 
which shows that the integrand is equivalent to 
\begin{equation*}
 \frac{C}{|\Re(s)|^{3/2}}e^{-(\alpha\log(\alpha) + (1-\alpha)\log(1-\alpha))\Re(s)}, \, |\Re(s)| \rightarrow \infty
 \end{equation*}
for some constant $C = C(\Im(s))$ depending only on $\Im(s)$.

\subsection{Behaviour for small indices: Cressie's result and its free analog} 
In this section, we investigate the limiting behaviour of $h_{\alpha}$ as $\alpha \rightarrow 0$. In the classical setting, the density tends to zero pointwisely therefore agrees somehow with Cressie's result we recall below. The latter implies that $a_{\alpha}(U)$ becomes degenerate when $\alpha$ approaches zero. In the free setting, the limiting density is an infinite positive measure whose image under the map $x \mapsto 1/x$ is the Haar measure on the compact interval $[0,1]$. Start with the series expansion of $h_{\alpha}$ is given by (A.6 p. 218 in \cite{HMS})
\begin{equation*}
h_{\alpha}(y) = (1-\alpha)\sum_{k=0}^{\infty} \frac{1}{\Gamma(1- \alpha - k\alpha )\Gamma(\alpha-k(1-\alpha))} \frac{(-1)^{k}}{k!} \frac{1}{y^{(k+1)(1-\alpha) +1}}
\end{equation*}
for $y > (1-\alpha)\alpha^{\alpha/(1-\alpha)}$. The mirror formula 
\begin{equation*}
\Gamma(\alpha-k(1-\alpha))\Gamma((1-\alpha)(k+1)) = \frac{\pi}{\sin((1-\alpha)(k+1)\pi)}
\end{equation*}
transforms $h_{\alpha}$ to 
\begin{equation*}
h_{\alpha}(y) = \frac{1-\alpha}{\pi} \sum_{k=0}^{\infty} \frac{\sin[(1-\alpha)(k+1)\pi] \Gamma[(1-\alpha)(k+1)]}{\Gamma(1- \alpha - k\alpha )} \frac{(-1)^k}{k!} \frac{1}{y^{(k+1)(1-\alpha) +1}}.
\end{equation*}
The defining term of the last series converges as $\alpha \rightarrow 0$ to 
\begin{equation*}
\frac{(-1)^k \sin((k+1)\pi) }{y^{k+2}} = 0
\end{equation*}
for any integer $k \geq 0$ while the support $((1-\alpha)\alpha^{\alpha/(1-\alpha)}, \infty)$ of $h_{\alpha}$ decreases to $[1,\infty)$ as $\alpha \rightarrow 0$.  Since 
\begin{equation*}
 \frac{\Gamma[(1-\alpha)(k+1)]}{\Gamma(1- \alpha - k\alpha )k!} \quad \rightarrow \quad 1
 \end{equation*}
as $\alpha \rightarrow 0$ and it defines a continuous function on the variable $\alpha$, then Lebesgue's convergence Theorem ensures the convergence of $h_{\alpha}(y)$, for fixed $y > 1$, to zero as the positive index $\alpha$ does. According to \cite{Cressie}, (see also \cite{Brockwell}), a stable variable with shape parameter $c=1$ and zero shift $\gamma = 0$ satisfies 
\begin{equation*}
(X_{\alpha})^{\alpha} \overset{d}{\rightarrow} \frac{1}{L}, \quad \alpha \rightarrow 0
\end{equation*}
so that $a_{\alpha}(U)$ converges in distribution to the Dirac measure $\delta_1$. For a positive free random variable, we start from the image of its density function under the map $x \mapsto x^{\alpha}$: 
\begin{equation*}
\frac{a_{1-\alpha}^{-1}(x^{1/\alpha})}{\pi x}\frac{\sin(\alpha a_{1-\alpha}^{-1}(x^{1/\alpha}))}{\alpha a_{1-\alpha}^{-1}(x^{1/\alpha})} \frac{\sin(a_{1-\alpha}^{-1}(x^{1/\alpha}))}{\sin((1-\alpha) a_{1-\alpha}^{-1}(x^{1/\alpha}))}, 
\end{equation*}
for $x > [a_{1-\alpha}(0)]^{\alpha} = \alpha^{\alpha}(1-\alpha)^{1-\alpha} \geq 1/2$. Then, we expand    
\begin{align*}
\frac{a_{1-\alpha}^{-1}(x^{1/\alpha})}{\pi} &=  2 H_{2,1}^{1,0}\left[\frac{1}{x^{2\alpha}}|_{(0,2)}^{(1,2\alpha),(1,2(1-\alpha))}\right] 
\\& = \sum_{k=0}^{\infty} \frac{1}{\Gamma(1- k\alpha )\Gamma(1- k(1-\alpha))} \frac{(-1)^{k}}{k!} \frac{1}{x^k}.
\end{align*}
With the help of the mirror formula, 
\begin{equation*}
\Gamma(1-k(1-\alpha))\Gamma((1-\alpha)k) = \frac{\pi}{\sin(k(1-\alpha)\pi)}, \quad k \geq 1, 
\end{equation*}  
one easily derives 
\begin{align*}
\lim_{\alpha \rightarrow 0} \frac{1}{\pi}a_{1-\alpha}^{-1}(x^{1/\alpha}) = 1 + \frac{1}{\pi}
\lim_{\alpha \rightarrow 0}\sum_{k=1}^{\infty} \frac{\sin(k(1-\alpha)\pi)\Gamma(k(1-\alpha))}{\Gamma(1- k\alpha)} \frac{(-1)^{k}}{k!} \frac{1}{x^k} = 1
 \end{align*}
for fixed $x > 1$. Since $\alpha^{\alpha}(1-\alpha)^{1-\alpha} \rightarrow 1$ as $\alpha \rightarrow 0$, then 
\begin{equation*}
\lim_{\alpha \rightarrow 0} \frac{a_{1-\alpha}^{-1}(x^{1/\alpha})}{\pi x}\frac{\sin(\alpha a_{1-\alpha}^{-1}(x^{1/\alpha}))}{\alpha a_{1-\alpha}^{-1}(x^{1/\alpha})} \frac{\sin(a_{1-\alpha}^{-1}(x^{1/\alpha}))}{\sin((1-\alpha) a_{1-\alpha}^{-1}(x^{1/\alpha}))} = \frac{1}{x} 
\end{equation*}
for any $x > 1$. Note that Scheff\'e's Lemma does not apply as in the classical setting, nevertheless the image of 
\begin{equation*}
\frac{dx}{x}{\bf 1}_{\{x > 1\}}
\end{equation*}
under the map $x \mapsto 1/x$ is the Haar measure on $[0,1]$, namely:
\begin{equation*}
\frac{dx}{x}{\bf 1}_{\{0 < x < 1\}}.
\end{equation*}

\section{On the function $a_{\alpha}$}
As stated in the introductory part, the occurence of $a_{\alpha}$ in both probabilistic settings is resumed in the integral representation of the density of $X_{\alpha}$ together with deformations of integration paths. Start with the inversion formula for the Fourier transform to see that the density of $X_{\alpha}$ reads (\cite{Ibr} p.49):
\begin{equation*}
\frac{1}{\pi} \Re \int_0^{\infty}e^{-itx} e^{-t^{\alpha}e^{-i\alpha\pi/2}} dt = \frac{1}{\pi} \Re \int_0^{\infty}e^{itx} e^{-(it)^{\alpha}} dt. 
\end{equation*}
The path of integration may be deformed into the positive half imaginary-line yielding (\cite{Zolo1})
\begin{align*}
\frac{1}{\pi} \Re \int_0^{\infty}e^{itx} e^{-(it)^{\alpha}} dt &= - \frac{1}{\pi}\Im  \int_0^{\infty}e^{-tx} e^{-e^{i\alpha\pi}t^{\alpha}} dt
\\& = \frac{1}{\pi}\Im  \int_0^{\infty}e^{-tx} e^{-e^{-i\alpha\pi}t^{\alpha}} dt.
\end{align*}
Performing the variable change $t=x^{1/(1-\alpha)}u$, the last integral transforms to 
\begin{equation*}
\frac{x^{1/(1-\alpha)}}{\pi} \Im \int_0^{\infty}e^{-x^{\alpha/(1-\alpha)} (u + e^{-i\alpha \pi}u^{\alpha})} du = \frac{x^{1/(1-\alpha)}}{\pi} \Im \int_0^{\infty}e^{-x^{\alpha/(1-\alpha)} (u + \phi_{1-\alpha}(u))} du
\end{equation*}
where $\phi_{1-\alpha}(u): = e^{-i\pi \alpha} u^{\alpha}$. Following \cite{Chernin} (see also Ch.II \cite{Ibr}), the path of integration may be taken as the curve 
\begin{align*}
\mathcal{C}_{\alpha} &:= \{u = re^{i\theta}, \theta \in [-\pi,0], \Im (u +\phi_{1-\alpha}(u)) = 0\}
\\& = \left\{u = re^{i\theta}, \theta \in [-\pi,0], r = \left[\frac{\sin(\alpha (\pi-\theta))}{\sin \theta}\right]^{1/(1-\alpha)}\right\}
\\& =  \left\{u = -re^{i\theta}, \theta \in [0,\pi], r = \left[\frac{\sin(\alpha \theta)}{\sin \theta}\right]^{1/(1-\alpha)}\right\}.
\end{align*}
Then for $u = -re^{i\theta} \in \mathcal{C}_{\alpha}$, 
\begin{align*}
\Re(u+\phi_{1-\alpha}(u)) &= -r\cos \theta + r^{\alpha} \cos(\alpha \theta) 
\\& = r^{\alpha}[-r^{1-\alpha}\cos \theta + \cos(\alpha\theta)] 
\\& = \left[\frac{\sin(\alpha \theta)}{\sin \theta}\right]^{\alpha/(1-\alpha)} \frac{\sin((1-\alpha)\theta)}{\sin(\theta)}
\end{align*}
and $a_{\alpha}(\theta)$ shows up. Now, consider $F_{\alpha}:= 1/G_{\alpha}$ where $G_{\alpha}$ is the Stieltjes  transform of the positive free stable distribution (\cite{Berco}). Then, $F_{\alpha}$ is a one-to-one correspondence from the upper half-plane $\mathbb{C}$ onto its image. Let $F_{\alpha}^{-1}$ be the reciprocal function of $F_{\alpha}$, then (\cite{Berco},\cite{Berco1}) 
\begin{equation*}
F_{\alpha}^{-1}(z) = \phi_{\alpha}(z) + z 
\end{equation*}
where $z$ belongs to the domain $\Im(z+ \phi_{\alpha}(z)) > 0$ and $\Im(z)$ is positive and sufficiently large. The function $\phi_{\alpha}$ is known as the Voiculescu transform of the positive free stable distribution. But since Stieltjes inversion formula involves the imaginary part of $G_{\alpha}$ when $z \in \mathbb{C}^+$ approaches the real line, one rather focuses on the domain  
\begin{equation*}
\Omega_{\alpha} := \{z, \, \Im(1/z+ \phi_{\alpha}(1/z) > 0\} 
\end{equation*}
for $z \in \mathbb{C}^-$ near the origin of the complex plane. This is a Jordan domain, that is delimited by a Jordan curve, whose boundary $\partial \Omega_{\alpha}$ is obviously
\begin{align*}
\partial \Omega_{\alpha} &= \left\{re^{i\theta}, \theta \in [-\pi,0],\, r = \left[\frac{\sin(\pi+ \theta)}{\sin ((1-\alpha)(\pi+\theta))}\right]^{1/\alpha}\right\}
\\& =  \left\{-re^{i\theta}, \theta \in [0,\pi],\, r = \left[\frac{\sin\theta}{\sin ((1-\alpha) \theta))}\right]^{1/\alpha}\right\}.
\end{align*}
This is the image of $\mathcal{C}_{1-\alpha}$ under the inversion map $z \mapsto 1/\overline{z}$. By Caratheodory's extension Theorem, $G_{\alpha}$ extends to a homeomorphism from the real line onto  $\partial \Omega_{\alpha}$ so that every $x$ lying in the support of the free positive stable variable is given by 
\begin{align*}
x &= G^{-1}_{\alpha}(z) ,\quad z \in  \partial \Omega_{\alpha}
\\& = \Re(1/z + \phi_{\alpha}(1/z)).
\end{align*}
The same computations performed in the classical setting lead to the function $a_{1-\alpha}$. \\
{\bf Acknowledgments}: the author would thanks W. Jedidi for making him aware of the paper \cite{PSS} and T. Simon for stimulating discussions held at Lille 1 university about stable random variables. He also gives a special thank to Professor J.J. Loeb for his remarks on the Mellin-Barnes integral defining the Fox H-function and to members of LAREMA at Angers university for their hospitality during his staying on June 2010.

 \end{document}